\documentclass[12pt]{amsart}

\usepackage{amssymb,latexsym}
\usepackage{verbatim,fullpage}

\input {cyracc.def}
\font\ququ=cmr10 scaled \magstep1
\font\tencyr=wncyr8 scaled \magstep1

\def\rus{\tencyr\cyracc}

\newcommand{\re}[1]{\textrm  (\ref{#1})}

\renewenvironment{proof}
{\noindent {\sl Proof.}\quad }{\hfill $\square$
\vskip1.1ex\noindent }

\newenvironment{proof*}
{\noindent {\sl Proof.}\quad }{\hfill $\square$}

\renewcommand{\theequation}{\thesection.\arabic{equation}}
\renewcommand{\thesubsubsection}{\theequation.\arabic{subsubsection}}

\catcode`\@=\active
\catcode`\@=11
\def\@eqnnum{\hbox to
.01pt{}\rlap{\hskip-\displaywidth(\mathbf{\theequation})}}

\catcode`\@=12

\newenvironment{s}[1]
{ \vskip1.2ex \refstepcounter{equation}
\noindent {\bf \theequation\quad #1.} \begin{sl}}{\end{sl}
\vskip1.1ex\noindent }

\newenvironment{rem}[1]
{ \vskip1.2ex \refstepcounter{equation}
\noindent {\bf \theequation\enspace {#1}.} }{ \vskip1.1ex\noindent }

\newcommand {\g}{{\frak g}}

\newcommand {\me}{{\frak m}}

\newcommand {\te}{{\frak t}}

\newcommand {\sln}{{\frak sl}_{n+1}}
\newcommand {\slv}{{\frak sl}(V)}

\newcommand {\esi}{\epsilon}
\newcommand {\ap}{\alpha}

\newcommand {\lb}{\lambda}

\newcommand {\vp}{\varphi}

\newcommand {\cc}{{\mathcal C}}
\newcommand {\ce}{{\mathcal E}}
\newcommand {\cf}{{\mathcal F}}
\newcommand {\ch}{{\mathcal H}}
\newcommand {\cj}{{\mathcal J}}
\newcommand {\N}{{\mathcal N}}
\newcommand {\cp}{{\mathcal P}}
\newcommand {\cR}{{\mathcal R}}

\newcommand {\cq}{{\mathcal Q}}
\newcommand {\VV}{{\Bbb V}}

\newcommand {\md}{/\!\!/}

\newcommand {\isom}{\stackrel{\sim}{\longrightarrow}}
\newcommand {\ad}{{\mathrm{ad\,}}}

\newcommand {\End}{{\mathrm{End\,}}}

\newcommand {\hot}{{\mathrm{ht}}}

\newcommand {\rk}{{\mathrm{rk\,}}}
\newcommand {\spe}{{\mathrm{Spec\,}}}

\newcommand {\tri}{{\frak sl}_2}
\newcommand {\GR}[2]{{\mathrm{{ #1}}}_{#2}}

\newcommand {\ov}{\overline}
\newcommand {\un}{\underline}

\newcommand{\PT}{{\mathcal P}_+}

\newcommand {\vno}[1]{\vskip#1 ex\noindent}

\newcommand {\beq}{\begin{equation}}
\newcommand {\eeq}{\end{equation}}

\newcommand{\evl}{\End\VV_\lb}
\newcommand{\etvl}{\End_T(\VV_\lb)}
\newcommand{\clg}{\mathcal C_\lb(\g)}
\newcommand{\clt}{\mathcal C_\lb(\te)}
\newcommand{\jlg}{J_\lb(\g)}
\newcommand{\dlq}{\mathcal D_\lb}
\newcommand{\fjlg}{\mathcal F(\jlg ; q)}
\newcommand{\mlm}{m_\lb^\mu}
\newcommand{\mlmq}{m_\lb^\mu(q)}
\newcommand{\vml}{V_\lb^\mu}
\newcommand{\vlb}{\VV_\lb}

\newfam\Bbbfam\newfam\eufam\newfam\Eufam%
\font\Bbbfont=msbm10 scaled 1200%
\font\olala=msam10 scaled 1200%
\font\frak=eufm10 scaled 1400%
\font\Bbbsmallfont=msbm8%
\textfont\Bbbfam=\Bbbfont\scriptfont\Bbbfam=\Bbbsmallfont
\font\euzw=eufm10 scaled 1200%
\font\euac=eufm7 scaled 1200%
\font\euacc=eufm7 scaled 1000%
\textfont\eufam=\euzw\scriptfont\eufam=\euac%
\scriptscriptfont\eufam=\euacc%
\font\Euzw=eufm10 scaled \magstep2%
\font\Euac=eufm7 scaled \magstep2%
\textfont\Eufam=\Euzw\scriptfont\Eufam=\Euac%
\def\frak{\fam\eufam}%
\def\Frak{\fam\Eufam}%
\def\Bbb{\fam\Bbbfam}%

\def\square{\hbox {\olala\char"03}}
\def\bbk{\hbox {\Bbbfont\char'174}}

\begin{document}
\setlength{\parskip}{2pt plus 4pt minus 0pt}
\hfill {\scriptsize March 26, 2002} \vskip1ex
\vskip1ex

\title[Endomorphism algebras and Dynkin polynomials]{
Weight multiplicity free representations,
$\Frak g$-endomorphism algebras,
and Dynkin polynomials 
}
\author
{\sc Dmitri I. Panyushev} 
\maketitle
\begin{center}
{\footnotesize
{\it Independent University of Moscow,
Bol'shoi Vlasevskii per. 11 \\
121002 Moscow, \quad Russia \\ e-mail}: {\tt panyush@mccme.ru }\\
}
\end{center}

\medskip
\smallskip

\section*{Introduction}
\vno{1}%
Throughout the paper, $G$ is a connected semisimple algebraic group defined
over an algebraically closed field $\bbk$ of characteristic zero,
and $\g$ is its Lie algebra.
\\[.6ex]
Recently, A.A.~Kirillov introduced an interesting class of associative
algebras connected with the adjoint representation of $G$ \cite{family}. In
our paper, such algebras are called $\g$-endomorphism algebras.
Each $\g$-endomorphism algebra is a module over the algebra of invariants
$\bbk[\g]^G$; furthermore, it is a direct sum of modules of covariants. Hence
it is a free graded finitely generated module over $\bbk[\g]^G$.
The aim of this paper is to show that
commutative $\g$-endomorphism algebras have intriguing connections with
representation theory, combinatorics, commutative algebra, and
equivariant cohomology.
\\[.6ex]
Let $\pi_\lb: G \to GL(\VV_\lb)$ be an irreducible representation, where
$\lb$ stands for the highest weight of $\VV_\lb$. Following Kirillov,
one can form an associative $\bbk$-algebra by taking the $G$-invariant
elements in the $G$-module $\End\VV_\lb\otimes \bbk[\g]$. That is, we set
\[
  \clg = (\End\VV_\lb\otimes \bbk[\g])^G \ .
\]
This algebra will be referred to as the {\it $\g$-endomorphism algebra\/} (of
type $\lb$). We do not use Kirillov's term `classical family algebra' for
$\clg$, nor we consider `quantum family algebras' in our paper.
It is proved in \cite{family} that $\clg$ is commutative if and only if
all weight spaces in $\VV_\lb$ are 1-dimensional. That paper also
contains a description of $\g$-endomorphisms algebras for simplest
representations of the classical Lie algebras.
In our paper we do not attempt to dwell upon consideration of particular cases,
but we rather try to investigate general properties of such algebras.
As a tool for studying $\g$-endomoprhism algebras, we use $\te$-endomorphism
algebras. Let $\te$ be a Cartan  subalgebra of $\g$ and $W$ the
corresponding Weyl group. Let $\End_T(\vlb)$ denote the set of $T$-equivariant
endomorphisms of $\vlb$, where $T$ is the maximal torus with Lie algebra $\te$.
Then
\[
 \clt= (\End_T(\vlb) \otimes \bbk[\te])^W
\]
is called the $\te$-{\it endomorphism algebra\/} of type $\lb$.
Both $\clt$ and $\clg$ are free graded $\bbk[\g]^G$-modules of the same rank,
and there exists an injective homomorphism of $\bbk[\g]^G$-modules
$\hat r_\lb: \clg\to \clt$.
Moreover, $\bbk[\g]^G$ is a subalgebra in both $\clt$ and  $\clg$, and
$\hat r_\lb$ is a monomorphism of $\bbk[\g]^G$-algebras. We show that
$\hat r_\lb$ becomes an isomorphism after inverting the discriminant
$D\in \bbk[\g]^G$.
\\[.6ex]
Most of our results concern the case in which $\clg$ is commutative, i.e.,
$\vlb$ is weight multiplicity free (={\sf wmf}). A considerable part of
{\sf wmf} representations consists of minuscule ones. We show that if
$\vlb$ is minuscule, then $\clg\simeq \bbk[\te]^{W_\lb}$. The proof relies
on a recent result of A.\,Broer concerning ``small'' $G$-modules \cite{bram}.
If $\vlb$ is {\sf wmf} but not minuscule, then both $\clt$ and $\clg$ have
zero-divisors. We prove that $\spe\clt$ is a disjoint union of affine spaces
of dimension $\rk\g$, while $\spe\clg$ is connected. But both
varieties have the same number of irreducible components, which is equal to
the number of dominant weights of $\vlb$. As a by-product, we obtain the
assertion that $\hat r_\lb$ is an isomorphism if and only if $\lb$ is
minuscule.  Since $\clg$ and $\clt$ are
graded $\bbk$-algebras, one may consider their Poincar\'e series. We explicitly
compute these series for any $\lb$.
\\[.6ex]
The principal result of the paper is that any commutative algebra $\clg$ is
Gorenstein. The proof goes as follows. Any set $f_1,\dots,f_l$
of algebraically independent homogeneous generators of $\bbk[\g]^G$ form a
system of parameters for $\clg$. Therefore $\clg$ is Gorenstein if and only
if $\cR^{(\lb)}:=\clg/\clg f_1+\dots +\clg f_l$ is. The finite-dimensional
$\bbk$-algebra $\cR^{(\lb)}$ is isomorphic with $(\evl)^A$, where $A\subset G$
is the connected centraliser of a regular nilpotent element.
Using this fact we prove
that the socle of $\cR^{(\lb)}$ is 1-dimensional, i.e.,
$\cR^{(\lb)}$ is Gorenstein. It is also shown that the Poincar\'e
polynomial of $\cR^{(\lb)}$ is equal to the Dynkin polynomial for $\vlb$.
The Dynkin polynomial is being defined for any $\vlb$. It can be
regarded as a $q$-analogue of $\dim\vlb$ that describes the distribution
of weight spaces with respect to some level function.
According to an old
result of E.B.\,Dynkin \cite{vereteno}, it is a symmetric unimodal polynomial
with integral coefficients. Later on, R.\,Stanley has observed that
Dynkin polynomials have rich combinatorial applications and
there is a multiplicative formula for them, see \cite{unimod}.
In our setting, the Dynkin polynomial of a {\sf wmf} $G$-module
$\vlb$  appears as the numerator of the Poincar\'e series of $\clg$.
It is natural to suspect that any reasonable finite-dimensional Gorenstein
$\bbk$-algebra is the cohomology algebra of a `good' variety. Following
this harmless idea, we construct for any {\sf wmf} representation of a
simple group $G$ a certain variety
$X_\lb\subset {\Bbb P}(\vlb)$. We conjecture that
$H^*(X_\lb)\simeq \cR^{(\lb)}$ and $H^*_{G_c}(X_\lb)\simeq\clg$, where
$G_c\subset G$ is a maximal compact subgroup.
If $\lb$ is a minuscule dominant weight, then $X_\lb$ is nothing but
$G/P_\lb$, a generalised flag variety. In this case, the conjecture
follows from the equality $\clg=\bbk[\te]^{W_\lb}$ (Theorem~\ref{main2})
and the well-known description of
$H^*_{G_c}(G/P_\lb)$. We also verify (a part of) the conjecture for
some non-minuscule weights.
\\[.6ex]
The paper is organised as follows. In Section~1, we collect necessary
information on modules of covariants. Section~2 is devoted to basic
properties of endomorphism algebras. We describe the structure of these
algebras in the commutative case, in particular, for the minuscule weights.
Dynkin polynomials and their applications are discussed in Section~3.
In Section~4, we give explicit formulas for the Poincar\'e series of
endomorphism algebras. The Gorenstein property is considered in Section~5.
Finally, in Section~6, we construct varieties $X_\lb\subset {\Bbb P}(\vlb)$
and discuss connections between $\g$-endomorphism algebras and
(equivariant) cohomology of $X_\lb$.

{\small {\bf Acknowledgements.}
I would like to thank Michel Brion for several useful remarks and
the suggestion to consider equivariant cohomology.
This research was supported in part by
RFBI Grant 01--01--00756.
}

\section{Generalities on modules of covariants}

\noindent
Fix a Borel subgroup $B\subset G$ and a maximal torus $T\subset B$.
We will always work
with roots, simple roots, positive roots, and dominant weights that are
determined by this
choice of the pair $(B,T)$. For instance, the roots of $B$ are positive and
a highest weight vector in some $G$-module is a $B$-eigenvector. More
specifically,
write $\cp$ for the $T$-weight lattice, and $\PT$ for the dominant
weights in $\cp$. Next, $\Delta$ (resp. $\Delta^+$) is the set of all
(resp. positive) roots, $\Pi$ is the set of simple roots,
and $\cq$ is the root lattice.
For a $G$-module $M$, let $M^\mu$ denote the $\mu$-weight space of $M$
($\mu\in\cp$). If $\lb\in\PT$, then $\VV_\lb$ stands for the simple
$G$-module with highest weight $\lb$. Set $\mlm=\dim \VV_\lb^\mu$.
The notation $\mu\dashv \VV_\lb$ means that $\mlm\ne 0$. For instance,
we have $0\dashv\VV_\lb$ if and only if $\lb\in\cq$.

Given $\lb\in\PT$, the space $(\VV_\lb\otimes\bbk[\g])^G$ is called the
{\it module of covariants\/} (of type $\lb$). We will write $\jlg$ for it.
Clearly, $\jlg$ is a module over $J_0(\g)=\bbk[\g]^G$,
and $\jlg\ne 0$ if and only if $\lb\in\cq$.
The elements of $\jlg$ can be identified with the $G$-equivariant morphisms
from $\g$ to $\VV_\lb$. More precisely, an element
$\sum v_i\otimes f_i\in\jlg$ defines the morphism
that takes $x\in\g$ to $\sum f_i(x)v_i\in\VV_\lb$. This interpretation of
$\jlg$ will freely be used in the sequel.
Since $\bbk[\g]^G$ is a graded algebra, each $\jlg$ is a graded module, too:
\[
  \jlg = \bigoplus_{n\ge 0}\jlg_n \ ,
\]
the component of grade $n$ being $(\VV_\lb\otimes\bbk[\g]_n)^G$.
The  {\it Poincar\'e series\/} of $\jlg$ is the formal power series
\[
  \fjlg=\sum_{n\ge 0} \dim\jlg_n q^n \in {\Bbb Z}[[q]] \ .
\]
More generally, $C$ being an arbitrary graded object, we write
$\mathcal F(C; q)$ for its Poincar\'e series.

The following fundamental result is due to B.~Kostant \cite{ko63}.

\begin{s}{Theorem}
$\jlg$ is a free graded $J_0(\g)$-module of rank $m_\lb^0$.
\end{s}%
Let $d_1,\dots,d_l$ be the degrees of basic invariants in $\bbk[\g]^G$,
where $l=\rk\g$.
It follows from the Theorem that $\fjlg$ is a rational function of the form
\[
   \fjlg =\frac{\sum_j q^{e_j(\lb)}}{\prod_{i=1}^l (1-q^{d_i})} \ .
\]
The numbers $\{e_j(\lb)\}$ $(1\le j\le m_\lb^0)$, which are merely the
degrees of a set of free homogeneous generators of $\jlg$, are called
the {\it generalised exponents\/} for $\VV_\lb$.
Another interpretation of generalised exponents is obtained as follows.
Let $f_1,\dots,f_l\in\bbk[\g]^G$ be a set of basic invariants, with
$\deg f_i=d_i$. It is a homogeneous system of parameters for $\jlg$.
Therefore
$\jlg=(f_1,\dots,f_l)\jlg\oplus \ch_\lb$, where $\ch_\lb$ is a graded
finite-dimensional $\bbk$-vector space such that $\dim\ch_\lb=m_\lb^0$.
Any homogeneous $\bbk$-basis for $\ch_\lb$ is also a basis for $\jlg$
as $\bbk[\g]^G$-module.
Let $\N\subset\g$ denote the nilpotent cone.
Since $\bbk[\g]^G/(f_1,\dots,f_l)\simeq \bbk[\N]$, we see that
\[
   \ch_\lb\simeq (\VV_\lb\otimes \bbk[\N])^G
     \quad\textrm{ and } \quad
   \mathcal F(\ch_\lb; q)=\sum_{i=1}^{m_\lb^0} q^{e_j(\lb)} \ .
\]
The polynomial $\mathcal F(\ch_\lb; q)$ has non-negative integral coefficients
and $\mathcal F(\ch_\lb; q)\vert_{q=1}= m_\lb^0$. It is a $q$-analogue
of $m_\lb^0$. A combinatorial formula for $\mathcal F(\ch_\lb; q)$ was found
by W.\,Hesselink \cite{wim} and D.\,Peterson, independently.
\\[.6ex]
Now we describe another approach to computing $\mathcal F(\ch_\lb; q)$,
which is due to
R.K.\,Brylinski. Let $e\in\N$ be a regular nilpotent element.
Then $G{\cdot}e$ is dense in $\N$ and $\bbk[\N]\simeq \bbk[G]^{G_e}$
\cite{ko63}. Hence $\ch_\lb\simeq (\vlb)^{G_e}$, in particular,
$\dim (\vlb)^{G_e}=m_\lb^0$.
Thus, the space $(\vlb)^{G_e}$ is equipped with a grading coming from the
above isomorphism. A direct description of this grading can be obtained in
terms
of ``jump polynomials''. Fix a principal $\tri$-triple $\{e,h,f\}$
such that $G_h=T$ and $e$ is a sum of root vectors corresponding
to the simple roots.
We have $G_e\subset B$ and $G_e\simeq Z(G)\times A$, where $Z(G)$
is the centre of $G$ and $A$ is a connected commutative unipotent group.
As $[h,e]=2e$, the space $(\vlb)^A$ is
$\ad h$-stable. It follows from the $\tri$-theory that $\ad h$-eigenvalues on
$(\vlb)^A$ are nonnegative. Moreover,
since $\lb\in\cq$, we see that $Z(G)$ acts trivially on $\vlb$,
$(\vlb)^{G_e}=(\vlb)^A$, and these eigenvalues are even.
Therefore it is convenient to consider $\tilde h=\frac{1}{2}h$ and its
eigenvalues. Set
\[
   (\vlb)^A_i=\{x\in (\vlb)^A\mid [\tilde h,x]=ix\} \ .
\]
and  $\cj_{\vlb}(q)=\sum_i \dim (\vlb)^A_i\, q^i$. This polynomial is called
the {\it jump polynomial\/} for $\vlb$.
\begin{s}{Theorem {\ququ \cite{jump},\cite{bram93}}} \label{jump}
For any $\lb\in \cq\cap\PT$, we have $\cj_{\vlb}(q)=\mathcal F(\ch_\lb; q)$.
\end{s}%
Let $\te$ denote Lie algebra of $T$ and
$W$ the Weyl group of $T$.
For any $\lb\in\cq$, the space $\VV_\lb^0$ is a $W$-module.
Therefore, one can form  the space $J_\lb(\te)=
(\VV_\lb^0\otimes\bbk[\te])^W$. It is a module over $J_0(\te)=\bbk[\te]^W$.
By Chevalley's theorem, the restriction homomorphism $\bbk[\g]\to\bbk[\te]$
induces an isomorphism
of  $J_0(\g)$ and $J_0(\te)$, so that this common algebra
will be denoted by $J$.
Since $W$ is a finite reflection group in $\te$,
we have $J_\lb(\te)$ is a free graded $J$-module of rank $m_\lb^0$.
Restricting a $G$-equivariant morphism $\g\to\VV_\lb$ to $\te\subset\g$
yields a $W$-equivariant morphism $\te\to \VV_\lb^0$. In other words, we
obtain a map $res_\lb:\jlg\to J_\lb(\te)$, which, in view of Chevalley's
theorem, is a homomorphism of $J$-modules. Since $G{\cdot}\te$ is dense in
$\g$,
the homomorphism $res_\lb$ is injective. It is not, however, always surjective.
The following elegant result is due to Bram Broer \cite[Theorem\,1]{bram}.

\begin{s}{Theorem} \label{bram}
Suppose $\lb\in \cq\cap\PT$. Then:\par
the homomorphism $res_\lb$ is onto $\Longleftrightarrow$
$m_\lb^{2\mu}=0$ for all $\mu\in\Delta$.
\end{s}%
In other words, $res_\lb$ is an isomorphism if and only if twice of a root is
not a weight for $\VV_\lb$. The $G$-modules satisfying the last condition
are said to be {\it small}.
\vskip1ex\noindent
Looking at the elements of $\jlg$ as $G$-equivariant morphisms
$\g\stackrel{\vp}{\longrightarrow} \VV_\lb$, one can consider the evaluation
map
\[
   \esi_x:  \jlg\to (\VV_\lb)^{G_x}
\]
for any $x\in\g$. Namely, set $\esi_x(\vp)=\vp(x)$. The following is a
particular case of a more general statement \cite[Theorem\,1]{ya},
which applies to arbitrary $G$-actions.

\begin{s}{Theorem} \label{moya}
Suppose $\ov{G{\cdot}x}$ is normal. Then
$\esi_x$ is onto.
\end{s}%

\section{$\g$-endomorphism and $\te$-endomorphism algebras}
\label{s-family}
\setcounter{equation}{0}
\noindent
Following A.A.\,Kirillov,
define the $\g$-{\it endomorphism algebra\/} of type $\lb$ by the formula
\[
  \clg = (\End\VV_\lb\otimes \bbk[\g])^G \ .
\]
It is immediate that $\clg$ is a $J$-module and an associative $\bbk$-algebra.
Since $\End\VV_\lb\simeq \VV_\lb\otimes\VV_\lb^\ast=\oplus c_\nu\VV_\nu$,
one sees that $\clg$ is a direct sum of modules of covariants (possibly
with multiplicities). Hence $\clg\simeq\oplus c_\nu J_\nu(\g)$.
It is important that all $\nu$ belong to $\cq$. Therefore
$\dim (\evl)^A=\dim (\evl)^0=\sum c_\nu m_\nu^0$.
Notice that $\clg$ is not only a $J$-module, but it also contains
$J$ as subalgebra, since $id_{\VV_\lb}\in\End\VV_\lb$. In particular,
$\clg$ is a $J$-algebra. Clearly, $\clg$ is a graded $\bbk$-algebra,
the component of grade $n$ being $(\evl\otimes\bbk[\g]_n)^G$.
\vno{1}%
The zero-weight space in the $G$-module
$\evl$ is the set of $T$-equivariant endomorphisms of $\VV_\lb$; i.e.,
we have $\etvl=(\evl)^0$.
Define the $\te$-{\it endomorphism algebra\/} of type $\lb$ by the formula
\[
  \clt= (\etvl\otimes \bbk[\te])^W \ .
\]
Using the above notation, one sees that
$\clt=\oplus c_\nu J_\nu(\te)$. It follows that, patching together the
homomorphisms $res_\nu$, one obtains the monomorphism of
$J$-algebras  $\hat r_\lb: \clg\to\clt$.
Thus, we have two associative $J$-algebras such that both are free graded
$J$-modules of the same rank, $\dim\etvl$.

\begin{s}{Lemma}   \label{discriminant}
Let $D\in J$ be the discriminant. Then $\clg_D$ and $\clt_D$ are isomorphic
as $J_D$-modules.
\end{s}\begin{proof}
Both $\clt$ and $\clg$ are built of modules of covariants. Therefore
the result stems from
the analogous statement for the modules of covariants, which has been proved
in \cite[Lemma\,1(iii)]{bram}; see also \cite[Prop.\,4]{ya} for another proof
in a more general context.
\end{proof}%
We will primarily be interested in commutative $\g$- and
$\te$-endomorphism algebras.
The following Proposition contains a criterion of commutativity. Part (i)
has been proved in \cite[Corollary \,1]{family}.
We give however a somewhat different proof for it, which has a potential of
being applied in more general situations, cf. \cite[Prop.\,3]{ya}.

\begin{s}{Proposition}
\begin{itemize}
\item[\sf (i)]
  $\clg$ is commutative $\Longleftrightarrow$ $\mlm=1$ for all $\mu\dashv
  \VV_\lb$;
\item[\sf (ii)]  $\clg$ is commutative $\Longleftrightarrow$ $\clt$ is
commutative.
\end{itemize}
\end{s}\begin{proof} (i)
``$\Leftarrow$''. \ Because $\etvl\simeq\oplus_{\mu}
\End(\vml)$,
the algebra $\clt$ is commutative whenever $\mlm=1$ for all $\mu\dashv\VV_\lb$.
Since $\hat r_\lb$ is injective, we are done.

``$\Rightarrow$''. \ Interpreting elements of $\clg$ as $G$-equivariant
morphisms, consider the evaluation map
\[
  \esi_x:\clg \to (\evl)^{G_x}   \quad (x\in\g) \ .
\]
Taking $x\in\te$ with $G_x=T$ and applying Theorem~\ref{moya}, we obtain
a surjective $\bbk$-algebra homomorphism $\clg\to\etvl=\oplus_\mu \End(\vml)$.
Thus, commutativity  of $\clg$ forces that of $\End(\vml)$ for
all $\mu\dashv\VV_\lb$. Whence $\mlm=1$.  \par
(ii) This readily follows from part (i) and from the fact that $\clg$ is a
subalgebra of $\clt$.
\end{proof}%
{\bf Definition.} A $G$-module $\VV_\lb$ is called {\it weight multiplicity
free} ({\sf wmf}), if $\mlm=1$ for all $\mu\dashv \VV_\lb$.
\vno{1}%
Although we do not need this directly, it is worth mentioning that
a classification of the {\sf wmf} $G$-modules with $G$ simple is
contained in \cite[4.6]{howe}.

\begin{s}{Lemma}  \label{wmf}
Suppose $\VV_\lb$ is {\sf wmf}. Then
$\End\VV_\lb$
is a multiplicity free $G$-module; that is, in the
decomposition $\End\VV_\lb=\oplus c_\nu\VV_\nu$ all nonzero
coefficients $c_\nu$ are equal to 1.
\end{s}\begin{proof}
An explicit formula for the multiplicities in tensor products is found in
\cite[Theorem\,2.1]{prv}. In particular, that formula says that
$c_\nu\le m_\lb^{\nu-\lb^*}\, (\le 1)$. Here $\lb^*$ is the highest weight
of the dual $G$-module $\VV_\lb^*$.
\end{proof}%
We regard roots and weights as elements of the $\Bbb Q$-vector space
$\cp\otimes_{\Bbb Z}{\Bbb Q}$ sitting in $\te^*$; next,
$(\ ,\ )$ is a fixed $W$-invariant bilinear form on $\te^*$ that is
positive-definite on $\cp\otimes_{\Bbb Z}{\Bbb Q}$. Then, as usual,
$\gamma^\vee=\frac{2\gamma}{(\gamma,\gamma)}$ for all $\gamma\in\Delta$.
\vno{.5}%
Of all {\sf wmf} $G$-modules, the minuscule ones occupy
a distinguished position. Recall that $\VV_\lb$ is called {\it minuscule}, if
$\lb$ is minuscule. Some properties of minuscule dominant
weights are presented in the following proposition, where either of two items
can be taken as a definition of a minuscule dominant weight.

\begin{s}{Proposition} \label{micro}
For $\lb\in\PT$, the following conditions are equivalent:
\begin{itemize}
\item[\sf (i)] \quad  If $\mu\dashv\VV_\lb$, then $\mu\in W\lb$;
\item[\sf (ii)] \quad $(\lb,\ap^\vee)\le 1$ for all $\ap\in\Delta^+$.
\end{itemize}
\end{s}%
It follows from Proposition~\ref{micro}(i)
that the corresponding simple $G$-module is {\sf wmf}, while
Proposition~\ref{micro}(ii) implies that if $\g$ is simple, then
a minuscule dominant weight is fundamental.

\begin{s}{Theorem} \label{main2}
If $\lb$ is minuscule, then $\clg\simeq \clt\simeq \bbk[\te]^{W_\lb}$.
\end{s}{\sl Proof.}\quad
1. Let us look again at the decomposition
\[
\VV_\lb\otimes\VV_\lb^*=\bigoplus_{\nu\in I}c_\nu \VV_\nu=
\VV_{\lb+\lb^*}\oplus(\bigoplus_{\nu < \lb+\lb^*}c_\nu\VV_\nu) \ .
\]
By Lemma~\ref{wmf}, all $c_\nu=1$. \par
2. Since $\lb$ is minuscule, each $\VV_\nu$ ($\nu\in I$) is
small in the sense of Broer. Indeed, we have $(\lb+\lb^*,\ap^\vee)\le
2$ for any $\ap\in\Delta^+$. Hence the same inequalities hold for an arbitrary
dominant weight $\mu\dashv (\VV_\lb\otimes\VV_\lb^*)$ in place of $\lb+\lb^*$.
On the other hand, if some $\VV_\nu$ were not small and hence $m_\nu^{2\gamma}
\ne 0$ for some $\gamma\in \Delta^+$, then we would obtain
$(2\gamma,\gamma^\vee)=4>2$. A contradiction!
\par
3. It follows from the previous part and Theorem~\ref{bram} that
$\hat r_\lb=\oplus_{\nu\in I} res_\nu :
\clg\to\clt$ is an isomorphism of $J$-algebras. Since all weight spaces in
$\VV_\lb$ are 1-dimensional and all weights are $W$-conjugate,
\[
 \etvl=\bigoplus_{\mu\in W\lb}\End(\vml) \simeq \bbk[W/W_\lb]
\]
as $W$-modules. Thus, \\
\hbox to \textwidth{\hfil
 $\clt =( \bbk[W/W_\lb]\otimes \bbk[\te])^W\simeq \bbk[\te]^{W_\lb}$.
 \hfil $\square$ }
\\[1.2ex]
{\bf Remark.} Another proof of this result, which does not appeal to
Broer's theorem, follows from the description of Poincar\'e series
given in Theorem~\ref{clg-poinc}.
\\
By a result of Steinberg, $W_\lb$ is a reflection group in $\te$. Hence
$\clg$ is a polynomial algebra, if $\VV_\lb$ is minuscule.
But for the other {\sf wmf} $G$-modules the situation  is not so good.
Since $\clg$ and $\clt$ are commutative $\bbk$-algebras in the {\sf wmf}-case,
one can consider the varieties $M_\lb(\te):=\spe\clt$ and
$M_\lb(\g):=\spe\clg$. The chain of algebras
$J\subset\clg\subset\clt$ yields a commutative diagram
$\begin{array}{ccc} M_\lb(\te)  &    & \\
                  \lefteqn{\tau\phantom{\mu}}\downarrow & \searrow & \\
                     M_\lb(\g)  & \to      & \g\md G\simeq \te /W
\end{array}$ . \ All maps here are finite flat morphisms.

\begin{s}{Theorem}  \label{2.6}
Let $\VV_\lb$ is a {\sf wmf} $G$-module. Then
\begin{itemize}
\item[{\sf (i)}] \ $\clt$ and $\clg$ are commutative Cohen-Macaulay
$\bbk$-algebras;
\item[{\sf (ii)}] \ $\clt$ and  $\clg$ are reduced;
\item[{\sf (iii)}] \ $M_\lb(\te)$ 
is a disjoint union of affine spaces of
dimension $l$. The connected (=\,irreducible) components of $M_\lb(\te)$
are parametrised
by the dominant weights in $\VV_\lb$;
\item[{\sf (iv)}] \  The morphism $\tau$ yields a one-to-one correspondence
between the irreducible components of $M_\lb(\te)$ and $M_\lb(\g)$.
The variety $M_\lb(\g)$ 
is connected.
\end{itemize}
\end{s}\begin{proof}
1. Cohen-Macaulayness follows, since both $\clt$ and $\clg$ are graded free
modules over the polynomial ring $J$. \par
2. Since $\clg$ is a subalgebra of $\clt$ via $\hat r_\lb$, it is enough to
prove that $\clt$ has no nilpotent elements. Let $\vp\in\clt$.
Regard $\vp$ as a $T$-equivariant morphism from $\te$ to
$\oplus_\mu \End(\vml)$.
If $\vp\ne 0$, then there is $\mu\dashv\VV_\lb$ such that
$\vp(t)(v_\mu)=c_\mu v_\mu$ for $t\in\te$, $0\ne v_\mu\in\vml$, and some
$c_\mu\in \bbk\setminus\{0\}$. Hence $\vp^n(t)(v_\mu)=(c_\mu)^n v_\mu\ne 0$.
\par
3. (cf. the proof of Theorem~\ref{main2}.)
Let $\mu_1,\dots,\mu_k$ be all dominant weights in $\VV_\lb$. Then
\[
  \bigoplus_\mu \End\vml=\bigoplus_{i=1}^k
\Bigl(\bigoplus_{w\in W}\End(\VV_\lb^{w\mu_i})\Bigr)=
\bigoplus_{i=1}^k\bbk[W/W_{\mu_i}]
\]
as $W$-modules. Hence
\[
  \clt = \Bigl((\bigoplus_{i=1}^k\bbk[W/W_{\mu_i}])\otimes\bbk[\te]\Bigr)^W=
\bigoplus_{i=1}^k (\bbk[W/W_{\mu_i}]\otimes \bbk[\te])^W=
\bigoplus_{i=1}^k \bbk[\te]^{W_{\mu_i}} \ .
\]
4. Because $\clt$ is a free $J$-module, none of the elements of $J$
becomes a zero-divisor in $\clt$ (and hence in $\clg$). Therefore, given any
$f\in J$, the principal open subset $M_\lb(\te)_f$ (resp. $M_\lb(\g)_f$) has
the same irreducible components as $M_\lb(\te)$ (resp. $M_\lb(\g)$).
Let us apply this to $f=D\in J$, the discriminant.
In the commutative case, one can consider $\clt_D$ not only as $J_D$-module,
which is a localisation of a $J$-module, but also as $\bbk$-algebra in its
own right. By Lemma~\ref{discriminant}, we then conclude that the
$\bbk$-{\it algebras\/} $\clt_D$ and $\clg_D$ are isomorphic.
\\
That $M_\lb(\g)$ is connected follows from the fact that $\clg$ is graded, and
the component of grade 0 is just $\bbk$. (Recall that  $\clg_n=(\evl\otimes
\bbk[\g]_n)^G$.)
\end{proof}%
We have shown that $\clt\simeq \oplus_{i=1}^k \bbk[\te]^{W_{\mu_i}}$. Since
$J\simeq \bbk[\te]^W$ and  $\bbk[\te]^W\subset \bbk[\te]^{W_{\mu_i}}$, it
is not hard to realise that $J$ embeds diagonally in the above sum.

\begin{s}{Corollary}
$\clg\simeq\clt$ if and only if $\lb$ is minuscule.
\end{s}\begin{proof}
One implication is proved in Theorem~\ref{main2}. Conversely,
if the algebras are isomorphic, then parts (iii) and (iv) show
$k=1$, that is, $\lb$ is minuscule.
\end{proof}%
In view of Theorem~\ref{bram}, the Corollary says that if $\lb$ is not
minuscule, then at least one simple $\g$-module occurring in $\VV_\lb
\otimes\VV_\lb^*$ is not small.

The discrepancy between $\clt$ and $\clg$ can be seen on the level of
Poincar\'e series. We will give below precise formulae for
$\mathcal F(\clt; q)$  and $\mathcal F(\clg; q)$.

\section{Dynkin polynomials} \label{ebd}
\setcounter{equation}{0}
\noindent
In 1950, E.B.\,Dynkin shown that to any simple $G$-module $\vlb$ one can
attach a symmetric unimodal polynomial \cite{vereteno}.
This polynomial represents the distribution of the weight spaces in
$\VV_\lb$ with respect to some level function. In Dynkin's note, the properties
of symmetricity and unimodality were also expressed by the words that
the weight system of $\vlb$ is ``spindle-like''. To prove this, Dynkin
introduced what is now called `a principal $\tri$-triple in $\g$'.
We say that the resulting polynomial is {\it the Dynkin polynomial\/}
(of type $\lb$).
In this section, we give some formulae for Dynkin polynomials and some
applications of them.
\vno{1}%
The lowest weight vector in $\vlb$ is $-\lb^*$. Recall that $\Pi$ is the set
of simple roots.
Given $\mu\dashv\vlb$, let us say that $\mu$ is on the $n$-th floor,
if $\mu-(-\lb^*)=\sum_{\ap\in\Pi}n_\ap \ap$ with $\sum n_\ap=n$.
Thus, the lowest weight is on the zero (ground) floor and the highest weight is
on the highest floor.

\noindent
{\bf Definition} (cf. \cite[p.222]{vereteno}).
Letting $(\vlb)_n=\displaystyle \sum_{\mu:\,floor(\mu)=n}\dim\vml$
and $a_n(\lb)=\dim(\vlb)_n$, define
the Dynkin polynomial $\dlq(q)$ to be $\sum_n a_n(\lb) q^n$.
(If we wish to explicitly indicate the dependence of this polynomial on $\g$, we
write $\dlq(\g)(q)$.
\\[.7ex]
To obtain a more formal presentation, consider the element
$\rho^\vee=\frac{1}{2}\sum_{\ap\in\Delta^+} \ap^\vee$. Since
$(\ap,\rho^\vee)=1$ for any $\ap\in\Pi$, we have $(\gamma,\rho^\vee)=
\hot\gamma$, the height of $\gamma$, for any $\gamma\in \Delta$.
If $\mu\in\cq$, then $(\mu,\rho^\vee)\in {\Bbb Z}$. Since
$\lb-(-\lb^*)\in\cq$ and $(\lb,\rho^\vee)=(\lb^*,\rho^\vee)$, we see
that $\hot(\mu):=(\mu,\rho^\vee)\in \frac{1}{2}{\Bbb Z}$ for an arbitrary
$\mu\in\cp$.
In view of these properties of $\rho^\vee$, the following is an obvious
reformulation of the previous definition:

\begin{equation}  \label{d-sum}
\dlq(q)=q^{(\lb,\rho^\vee)}\sum_{\mu\dashv \vlb}\mlm\, q^{(\mu,\rho^\vee)}
=\sum_{\mu\dashv \vlb}\mlm\, q^{\hot(\lb+\mu)}  \ .
\end{equation}
Thus, $\rho^\vee$ defines a grading of $\vlb$, and $\dlq(q)$ is nothing
but the shifted Poincar\'e polynomial of this grading.

\begin{s}{Theorem {\ququ \cite[Theorem\,4]{vereteno}}}
For any $\lb\in\PT$, the polynomial $\dlq$ is symmetric
(i.e. $a_i(\lb)=a_{m-i}(\lb)$, $m=\deg\dlq$) and unimodal
(i.e. $a_0(\lb)\le a_1(\lb)\le\dots \le a_{[m/2]}(\lb)$).
\end{s}%
{\sl Idea of proof.}\quad 
Having identified $\te$ and $\te^*$, we see that $2\rho^\vee$
becomes $h$, the semisimple element of our fixed principal $\tri$-triple.
Therefore, $(\vlb)_n=\{v\in\vlb \mid h{\cdot}v=2(n-\hot(\lb))v\}$.
Hence, up to a shift of degree, $\dlq(q^2)$ gives the character of $\vlb$ as
$\tri$-module. From \re{d-sum} one also sees that $\deg\dlq=2(\rho^\vee,\lb)
=2\,\hot(\lb)$.

The next proposition provides a multiplicative formula for the Dynkin
polynomials.
Apparently, this formula was first proved, in a slightly
different form, in \cite{unimod}. It was R.\,Stanley
who realised that Dynkin's result has numerous combinatorial
applications. Of course, the proof exploits Weyl's character formula.

\begin{s}{Proposition}  \label{d-mult}
$ \displaystyle \dlq(q)=\prod_{\ap\in\Delta^+}
  \frac{1-q^{(\rho+\lb,\ap^\vee)}}{1-q^{(\rho,\ap^\vee) }}$, where
  $\rho=\frac{1}{2}\sum_{\ap\in\Delta^+}\ap$.
\end{s}%
{\bf Remark.} A similar formula also appears in \cite[Lemma\,2.5]{stembridge}
(and probably in many other places) as a $q$-analogue of Weyl's
dimension formula (or a specialization of Weyl's
character formula). However, Stembridge did not make the degree shift
$q^{(\rho^\vee,\lb)}$ and did not mention a connection of this $q$-analogue
with Dynkin's results.

\begin{rem}{Example}  \label{primer}
Suppose $\g=\sln$, and let $\vp_i$ be the $i$-th fundamental weight of it. Then
a direct calculation based on \re{d-mult} gives:
\[
  \mathcal D_{m\vp_1}(\sln)(q)=
\frac{(1-q^{m+1})\dots (1-q^{m+n})}{(1-q)\dots (1-q^{n})}=:
\genfrac{[}{]}{0pt}{}{m+n}{n}=\genfrac{[}{]}{0pt}{}{m+n}{m} \ ,
\]
\begin{center}
  $\mathcal D_{\vp_m}(\sln)(q)=
\displaystyle \genfrac{[}{]}{0pt}{}{n+1}{m}$ \ .
\end{center}
\end{rem}%
As a consequence of this example, one can deduce the following (well-known)
assertions:
\begin{s}{Proposition}  \label{dlq}
Suppose $\dim V=2$. Then there are two isomorphisms of $\slv$-modules:
\[   S^n(S^m(V))\simeq S^m(S^n(V))\quad \textrm{ (Hermite's reciprocity)}
\]
\hbox to \textwidth{\hfil   $\wedge^m (S^{n+m-1}(V))=S^m(S^n(V))$. \hfil}
\vskip-1ex
\end{s}\begin{proof}
The first formula follows from the equality
$\mathcal D_{m\vp_1}(\sln)=\mathcal D_{n\vp_1}({\frak sl}_{m+1})$ and the
fact that $\VV_{m\vp_1}(\sln)\vert_{\tri}=S^m(S^n(V))$, where $\tri=\slv$ is a
principal $\tri$ in $\sln$.
Similarly, the second formula follows from the equality
$\mathcal D_{m\vp_1}(\sln)=\mathcal D_{\vp_m}({\frak sl}_{m+n})$.
\end{proof}%
Now we show that Dynkin polynomials arise in connection with
$\g$-endomorphism algebras for the minuscule dominant weights.
\\
For any $\lb\in\PT$, we have $W_\lb$ is a parabolic subgroup of $W$.
Let $d_i(W_\lb)$, $1\le i\le l$,  be the degrees of basic invariants in
$\bbk[\te]^{W_\lb}$. In particular, $d_i=d_i(W_0)=d_i(W)$.
Note that if $\lb\ne 0$, then some $d_i(W_\lb)$ are equal to 1.
Let $n: W\to {\Bbb N}\cup \{0\}$ be the length function with respect to
the set of simple reflections. Set $t_\lb(q)=
\sum _{w\in W_\lb} q^{n(w)}$. It is well known (see e.g. \cite[3.15]{hump})
that
\begin{equation}  \label{tlq}
 t_\lb(q)=\prod_{i=1}^l \frac{ 1-q^{d_i(W_\lb)} }{1-q} \ .
\end{equation}
These polynomials will frequently appear in the following exposition.

\begin{s}{Proposition}  \label{dlq-minus}
Suppose $\lb\in\PT$ is minuscule. Then \\[.6ex]
\hbox to \textwidth{\hfil
$\displaystyle \dlq(q)=\frac{t_0(q)}{t_\lb(q)}= \prod_{i=1}^l
\frac{ 1-q^{d_i} }{1-q^{d_i(W_\lb)}}$ .
\hfil}
\vskip-1ex
\end{s}\begin{proof}
Set $\Delta(\lb)=\{\ap\in\Delta\mid (\ap,\lb)=0\}$.
Using Propositions~\ref{d-mult} and \ref{micro}(ii), we obtain
\[
 \dlq(q)=\prod_{\ap\in\Delta^+\setminus \Delta(\lb)^+}
  \frac{1-q^{(\rho,\ap^\vee)+1}}{1-q^{(\rho,\ap^\vee) }}=
\prod_{\ap\in\Delta^+}
  \frac{1-q^{(\rho,\ap^\vee)+1}}{1-q^{(\rho,\ap^\vee) }}\cdot
\prod_{\ap\in\Delta(\lb)^+}
  \frac{1-q^{(\rho,\ap^\vee) }}{1-q^{(\rho,\ap^\vee)+1}} \ .
\]
Note that $(\rho,\ap^\vee)=\hot\ap^\vee$, the height of $\ap^\vee$ in the
dual root system $\Delta^\vee$. By a result of Kostant (see \cite[3.20]{hump}),
\[
   t_0(q)=\prod_{\ap\in\Delta^+}
\frac{ 1-q^{\hot(\ap)+1} }{ 1-q^{\hot(\ap)} } \ .
\]
Applying this formula to $\Delta^\vee$ and $\Delta(\lb)^\vee$
and substituting in the previous expression for $\dlq$, we
complete the proof.
\end{proof}%
Let us look at the Poincar\'e series of $\clg$, where $\lb$ is minuscule.
Since $\clg=\bbk[\te]^{W_\lb}$ (Theorem~\ref{main2}) and
$\mathcal F(J;q)=\prod_{i=1}^l (1/(1-q^{d_i}))$, we deduce from
Proposition~\ref{dlq-minus} that
\begin{equation} \label{poinc-minus}
 \mathcal F(\clg;q) 
=\dlq(q)\cdot \mathcal F(J; q) \ .
\end{equation}
It will be shown in Section~\ref{gor} that
this relation holds for all
{\sf wmf} $G$-modules.

\section{The Poincar\'e series of endomorphism algebras}
\setcounter{equation}{0}
\noindent
In this section we find explicit formulas for the
Poincar\'e series $\cf(\clg; q)$ and $\cf(\clt; q)$ with arbitrary
$\lb\in\PT$.
\\[.6ex]
Let $\mlmq$ be Lusztig's $q$-analogue of weight multiplicity. It is a certain
polynomial in $q$, with integral coefficients,  such that $\mlm(1)=\mlm$.
Defining of these $q$-analogues consists of two steps. First, one
defines a $q$-analogue of Kostant's partition function ${\mathcal P}$ by
\[
 \prod_{\ap\in\Delta^+}\frac{1}{1-e^\ap q}=\sum_{\nu\in\cq}
{\mathcal P}_q(\nu)e^\nu \ .
\]
Then we set $\mlmq=\sum_{w\in W} \det(w){\mathcal P}_q(w(\lb+\rho)-\mu-\rho)$.
Proofs of the next properties of polynomials $\mlmq$ can be found in
\cite{bram93} and \cite{ranee87}.
\begin{itemize}
\item If both $\lb,\mu$ are dominant, then
   \begin{itemize}
     \item the coefficients of $\mlmq$ are nonnegative;
     \item $\mlmq\ne 0$ $\Leftrightarrow$ $\mu\dashv\vlb$;
     \item $\deg\mlmq=(\lb-\mu,\rho^\vee)=\hot(\lb-\mu)$;
   \end{itemize}
\item  If $\lb\in\cq$, then $m_\lb^0(q)$ is the numerator of the Poincar\'e
series for the module of covariants
$\jlg$; that is, in the notation of Section~1, we have
$\cf(\ch_\lb;q)=\cj_{\vlb}(q)=m_\lb^0(q)$.
\end{itemize}
\begin{s}{Theorem}  \label{clg-poinc}
Let $\lb$ be an arbitrary dominant weight.
\begin{itemize}
\item[\sf (i)]  $\cf(\clg; q)=\displaystyle
\sum_{\nu\in\PT} (m_\lb^\nu(q))^2\,\frac{t_0(q)}{t_\nu(q)}\left/
\prod_{i=1}^l (1-q^{d_i})\right.=
\sum_{\nu\in\PT} (m_\lb^\nu(q))^2\,\frac{t_0(q)}{t_\nu(q)}\cdot \cf(J;q)$;
\item[\sf (ii)] $\cf(\clt; q)=\displaystyle
\sum_{\nu\dashv\vlb,\,\nu\in\PT}\frac{t_0(q)}{t_\nu(q)}\left/
\prod_{i=1}^l (1-q^{d_i})\right.=
\sum_{\nu\dashv\vlb,\,\nu\in\PT}\frac{t_0(q)}{t_\nu(q)}\cdot \cf(J;q)$.
\end{itemize}
\end{s}\begin{proof}
(i) Recall that $J=\bbk[f_1,\dots,f_l]$ and $\deg f_i=d_i$. Since
$f_1,\dots,f_l$ is a homogeneous system of parameters for
$\clg$ as $J$-module, we have
\[
\clg/\clg f_1+\dots \clg f_l\simeq (\evl\otimes\bbk[\N])^G\simeq (\evl)^A
\]
(cf. Section~1).
Hence $\cf(\clg; q)=\cf((\evl)^A;q)/\prod_{i=1}^l (1-q^{d_i})$. Using the
decomposition $\evl=\oplus_{\nu\in I} c_\nu\VV_\nu$ and
Theorem~\ref{jump}, we see
that the numerator is just the jump polynomial corresponding to the
(reducible) $G$-module $\evl$.
That is, $\cf((\evl)^A;q)=\sum c_\nu m_\nu^0(q)$.
It is remarkable however that, for the $G$-modules
of the form $\vlb\otimes \VV_\mu^*$, there is a formula for the jump polynomial
that does not appeal to the explicit decomposition of this tensor product.
Namely, by
\cite[Corollary\,2.4]{ranee}, we have
\begin{equation}  \label{ranee-jump}
\cj_{\vlb\otimes\VV_\mu^*}(q) =
\sum_{\nu\in\PT} m_\lb^\nu(q)m_\mu^\nu(q)\,\frac{t_0(q)}{t_\nu(q)} \ .
\end{equation}
Letting $\mu=\lb$, one obtains the required formula.
\par
(ii)  Recall from Section~2 that $\clt\simeq \oplus_\nu
\bbk[\te]^{W_{\nu}}$, where $\nu$ ranges over
all dominant weights of $\vlb$. Making use of Eq.~\re{tlq}, we obtain
\[
\cf(\bbk[\te]^{W_{\nu}};q)=1/\prod_{i=1}^l (1-q^{d_i(W_\nu)})=
\frac{t_0(q)}{t_\nu(q)}\left/\prod_{i=1}^l (1-q^{d_i})\right. \ ,
\]
which completes the proof.
\end{proof}%
In case $\lb$ is minuscule, Theorem~\ref{clg-poinc} shows that
\[
 \cf(\clg;q)=\cf(\clt;q)=\frac{t_0(q)}{t_\lb(q)}\cdot\cf (J;q) \ ,
\]
since $m_\lb^\lb(q)=1$. Thus, we recover in this way Eq.~\re{poinc-minus}
and the fact that $\clg=\clt$.
>From the last equality, we deduce the
following claim.  (This yields another proof for a part Broer's
results.)
\begin{s}{Corollary} Suppose $\lb$ is minuscule, and let $\VV_\nu$ be any
irreducible constituent of $\vlb\otimes\VV_\lb^*$. Then the restriction
homomorphism $res_\nu :J_\nu(\g) \to J_\nu(\te)$ is onto.
\end{s}%
Let $\cf_\lb(q)$ denote the right-hand side of Eq.~\re{ranee-jump} with
$\lb=\mu$,
i.e., the jump polynomial for the $G$-module $\evl$.

\begin{s}{Lemma}   \label{deg-flq}
\begin{itemize}
\item[\sf (i)] $\deg\cf_\lb(q)= 2\,\hot(\lb)$;
\item[\sf (ii)] $\cf_\lb(1)=\sum_{\mu\dashv\vlb} (\mlm)^2$;
\item[\sf (iii)]
If $\vlb$ is\/ {\sf wmf}, then $\displaystyle
\cf_\lb(q)= \sum_{\mu\dashv\vlb,\,\mu\in\PT} q^{2\hot(\lb-\mu)}\cdot
\frac{t_0(q)}{t_\mu(q)}$\, .
\end{itemize}
\end{s}\begin{proof}
1. $\cf_\lb(q)$ is a sum (with multiplicities) of the jump polynomials
for the irreducible constituents of $\evl$.
The jump polynomial for $\VV_{\lb+\lb^*}$ is of degree
$\hot(\lb+\lb^*)=2\hot(\lb)$. For all other simple $G$-submodules in
$\evl$, the height of the highest weight is strictly less.
\par
2. $\cf_\lb(1)=\dim(\evl)^A=\dim(\evl)^T=\sum_{\mu}\dim\End(\vml)$.
\par
3. In view of \ref{clg-poinc}(i), it suffices to prove that
$\mlmq=q^{\hot(\lb-\mu)}$.
Because $\mlm=1$ and the coefficients of $\mlmq$ are non-negative integers,
$\mlmq=q^a$. Since $\deg\mlmq=\hot(\lb-\mu)$, we are done.
\end{proof}%
The last expression demonstrates an advantage
of using  Eq.~\re{ranee-jump} in the {\sf wmf} case. We obtain a
closed formula for the jump polynomial of a reducible representation
that requires no bulky computations.


\section{The Gorenstein property for the commutative $\g$-endomorphism
algebras}  \label{gor}
\setcounter{equation}{0}
\noindent
The goal of this section is to prove that if $\vlb$ is {\sf wmf}, then
$\clg$ is a Gorenstein algebra. We also give another expression
for the Poincare series of $\clg$,
which includes the Dynkin polynomial of type $\lb$.
\\[.6ex]
To begin with, we recall some facts on graded Gorenstein algebras. A nice
exposition of relevant material is found in \cite{combin}.
Let $\cc=\oplus_{n\ge 0} \cc_n$ be a graded Cohen-Macaulay $\bbk$-algebra
with $\cc_0=\bbk$.
Suppose the Krull dimension of $\cc$ is $n$ and let $f_1,\dots,f_n$ be
a homogeneous system of parameters. Then $\ov{\cc}=\cc/(f_1,\dots,f_n)$
is a graded Artinian Cohen-Macaulay $\bbk$-algebra. We have
$\ov{\cc}=\oplus_{i=0}^d \ov{\cc}_i$ for some $d$, and
$\me=\oplus_{i=1}^d \ov{\cc}_i$ is the unique maximal ideal in $\ov{\cc}$.
The annihilator of $\me$ in $\ov{\cc}$ is called the {\it socle\/} of
$\ov{\cc}$:
\[
 \mathrm{soc}(\ov{\cc})=\{c\in \ov{\cc}\mid c{\cdot}\me=0\} \ .
\]
Then the following is true:
\\[.6ex]\indent
$\bullet$\quad $\cc$ is Gorenstein if and only if $\ov{\cc}$ is;
\par
$\bullet$\quad $\ov{\cc}$ is Gorenstein if and only if
$\mathrm{soc}(\ov{\cc})$ is 1-dimensional.
\\[.6ex]
At the rest of the section, $\vlb$ is a {\sf wmf} $G$-module, and hence
$\clg$ is commutative.
Consider the $\bbk$-algebra $(\evl)^A=:\cR^{(\lb)}$. Being a homomorphic image
of $\clg$, it is commutative as well. We also have $\dim\cR^{(\lb)}=\dim\vlb$.

\begin{s}{Proposition}  \label{privet}
Let $v_{-\lb^*}$ be a lowest weight vector in $\vlb$. Then
$\cR^{(\lb)}(v_{-\lb^*})=\vlb$.
\end{s}\begin{proof*}
Set $N=\cR^{(\lb)}(-v_{\lb^*})$, and consider $N^\perp$, the annihilator of $N$ in the
dual space $\VV_{\lb^*}$.  Let $v_{\lb^*}$ be a highest weight vector in $\VV_{\lb^*}$.
By definition, $N$ contains $v_{-\lb^*}$ and therefore $v_{\lb^*}\not\in
N^\perp$.
Since $A\subset G\subset\evl$ and $A$ is commutative,  $N$ is an $A$-module
and hence $N^\perp$ is an $A$-module, too.
Assume that $N^\perp\ne 0$. Since $A$ is unipotent, $N^\perp$ must contain
a non-trivial $A$-fixed vector. By a result of Graham (see \cite[1.6]{graham}),
$\dim(\VV_{\lb^*})^A=1$ in the {\sf wmf} case. As $A\subset B$, we conclude
that $(\VV_{\lb^*})^A= \bbk v_{\lb^*}$. Thus, $N^\perp\cap (\VV_{\lb^*})^A=0$.
This contradiction shows that $N^\perp=0$.
\end{proof*}%
\begin{s}{Proposition} \label{privet2}
$\cf_\lb(q)=\dlq(q)$ if and only if\/ $\vlb$ is {\sf wmf}.
In particular, $\cf_\lb(q)$ is symmetric and unimodal in the
{\sf wmf}-case.
\end{s}\begin{proof}
1. ``$\Rightarrow$''. Suppose $\cf_\lb(q)=\dlq(q)$. Then
\[
\sum_{\mu\dashv\vlb} (\mlm)^2=\cf_\lb(1)=\dlq(1)=\dim\vlb=
\sum_{\mu\dashv\vlb} \mlm \ .
\]
Whence $\vlb$ is {\sf wmf}. \par
2. ``$\Leftarrow$''.
The polynomial $\cf_\lb(q)$ is determined via the $\frac{1}{2}h$-grading
in $\cR^{(\lb)}=(\evl)^A$, whereas $\dlq(q)$ is determined via the shifted
`$\rho^\vee$-grading' in $\vlb$. Identifying $\te$ and $\te^*$, we obtain
$h=2\rho^\vee$. Obviously,  the bijective linear map
$\cR^{(\lb)}\to \vlb$, $x\mapsto x(v_{-\lb^*})$, respects both
gradings and has degree zero. Hence
$\cR^{(\lb)}_i\isom (\vlb)_i$.
\end{proof}%
{\bf Remark.} If $\vlb$ is not {\sf wmf}, then $\cf_\lb(q)$ can be neither
symmetric, nor unimodal. For instance, if $\g={\frak sp}_6$ and $\lb=\vp_2$,
the second fundamental weight, then $\cf_{\vp_2}(q)=
1+q+2q^2+2q^3+3q^4+2q^5+3q^6+q^7+q^8$.
\\[1ex]
By a theorem of Stanley \cite[12.7]{combin},
if $\cc$ is a Cohen-Macaulay \un{domain},
then the symmetricity of the Poincar\'e polynomial of
$\cc/(f_1,\dots,f_n)$ implies the Gorenstein property for $\cc$.
In our situation, $\clg$ is not a domain unless $\lb$ is minuscule,
see Theorem~\ref{2.6}.
So, we still cannot conclude that $\clg$ is always Gorenstein.

\begin{s}{Proposition}  \label{privet3}
The socle of $\cR^{(\lb)}$   
is 1-dimensional.
\end{s}\begin{proof}
Write $\cR$ for $\cR^{(\lb)}$ in this proof.
Recall from Section~\ref{ebd} that $\vlb$ is a disjoint union of `floors',
the weight space $\vlb^{-\lb^*}$ being the zero floor and $\vlb^\lb$ being
the $2\,\hot(\lb)$-th floor.
We know that $\cR$ is 
commutative 
and $\cR=\displaystyle \bigoplus_{i=0}^{2\hot(\lb)}\cR_i$.
By Propositions~\ref{privet} and \ref{privet2}, we have
$\cR_i(v_{-\lb^*})=(\vlb)_i$, the $i$-th floor in $\vlb$, and
$\dim \cR_{2\,\hot(\lb)}=1$. Clearly, $\cR_{2\,\hot(\lb)}$
takes $\vlb^{-\lb^*}$ to $\vlb^\lb$, and $\cR_{2\,\hot(\lb)}(\vml)=0$
for $\mu\ne -\lb^*$. We wish to show that
soc$(\cR)=\cR_{2\,\hot(\lb)}$.
Suppose $x\in \cR_i$;
then $xv_{-\lb^*}\in (\vlb)_i$. Let $\langle\ ,\ \rangle$ denote
the natural pairing of $\vlb$ and $\VV_{\lb^*}$. Using the $h$-invariance,
we see that $\langle (\vlb)_i,(\VV_{\lb^*})_j\rangle=0$ unless $i+j=2\,\hot(\lb)$.
Hence
there exists $\xi\in (\VV_{\lb^*})_j$, where $j={2\,\hot(\lb)-i}$, such that
$\langle\xi,  xv_{-\lb^*}\rangle\ne 0$. Applying Proposition~\ref{privet}
to $\VV_{\lb^*}$, we obtain $\xi=y(v_{-\lb})$ for some $y\in \cR_j$.
Here $v_{-\lb}$ is a lowest weight vector in $\VV_{\lb^*}$, and
$\cR$ is identified with $(\End\VV_{\lb^*})^A$. Hence
$\langle v_{-\lb}, yx(v_{-\lb^*})\rangle\ne 0$.
Thus, for any $x\in\cR_i$ there exists $y\in \cR_{2\,\hot(\lb)-i}$ such
that $xy\ne 0$, which is exactly what we need.
\end{proof}%
Combining all previous results of this section, we obtain

\begin{s}{Theorem}  \label{main5}
Let $\vlb$ be a\/ {\sf wmf} $G$-module.  Then
\begin{itemize}
\item[\sf (i)] \quad  The map $(\evl)^A \to \vlb$, $x\mapsto x(v_{-\lb^*})$,
is bijective;
\item[\sf (ii)] \quad $(\evl)^A$ is an Artinian Gorenstein $\bbk$-algebra;
\item[\sf (iii)] \quad $\dlq(q)=\cf_\lb(q)$;
\item[\sf (iv)] \quad $\clg$ is a Gorenstein $\bbk$-algebra.
\end{itemize}
\end{s}%
It is worth mentioning the following property of {\sf wmf} $G$-modules
whose proof also uses Graham's result, cf. Proposition~\ref{privet}.

\begin{s}{Proposition}
For any $n\in \Bbb N$, the vector $e^n(v_{-\lb^*})$ has nonzero projections
to all weight spaces in $(\vlb)_n$.
\end{s}\begin{proof}
Let $\ce$ be the $\bbk$-subalgebra of $\cR^{(\lb)}$ generated by $e$.
Then $M:= (\ce(v_{-\lb^*}))^\perp$ is an $\ce$-stable subspace of
$\VV_{\lb^*}$, and, obviously, $v_{\lb^*}\not\in M$. Assume that
$\ce(v_{-\lb^*})$ has the zero projection to some weight space, say
$\vlb^{-\nu}$. Then $M$ contains the weight space $\VV_{\lb^*}^\nu$. Let
$v_\nu$ be a nonzero vector in $\VV_{\lb^*}^\nu$ and let $k$ be
the maximal integer such that $e^k(v_\nu)\ne 0$. By \cite[2.6]{jump},
$e^k(v_\nu)\in (\VV_{\lb^*})^A$ and by \cite[1.6]{graham},
$(\VV_{\lb^*})^A=\bbk v_{\lb^*}$. As $e^k(v_\nu)\in M$, we obtain a
contradiction.
\end{proof}%
%
All {\sf wmf} representation of simple
algebraic groups are found by R.\,Howe in \cite[4.6]{howe}.
This information is included in the first two columns of Table~1;
the last column gives the corresponding Dynkin polynomial.
We write $\vp_i$ for the $i$-th fundamental weight of $G$, according to
the numbering of \cite{vion}.

\begin{table}[htb]
\begin{center}
\caption{ 
Dynkin polynomials for the weight
multiplicity free representations}
\begin{tabular}{cccc}
   Type        & $\lb$     & minuscule &  $\dlq(q)$  \\ \hline\hline
\smash{\raisebox{-2.5ex}{$\GR{A}{n}$}}   & $\vp_i$   & yes  & see Example~\ref{primer}
\\ 
          &  $\begin{array}{c} m\vp_1, m\vp_n \\
             {\scriptstyle (m\ge 2)} \end{array}$  & no & see Example~\ref{primer}\\ \hline
\smash{\raisebox{-2ex}{$\GR{B}{n}$}} & $\vp_n$ & yes & $(1+q)(1+q^2)\dots (1+q^n)$
\\ 
            & $\vp_1$ & no  & $1+q+\dots +q^{2n}$ \\ \hline
\smash{\raisebox{-2ex}{$\GR{C}{n}$}} & $\vp_1$ & yes  & $1+q+\dots +q^{2n-1}$
\\ 
            & $\vp_3\ (n=3)$ & no & $1+q+q^2+2(q^3+\dots +q^6)+q^7+q^8+q^9$
                      \\ \hline
\smash{\raisebox{-2ex}{$\GR{D}{n}$}} & $\vp_1$ & yes & $(1+q^{n-1})(1+q+\dots +q^{n-1})$
\\ 
    & $\vp_{n-1},\vp_n$ & yes & $(1+q)(1+q^2)\dots (1+q^{n-1})$ \\ \hline
$\GR{E}{6}$  & $\vp_1$  & yes & $(1+q^4+q^8)(1+q+\dots +q^8)$ \\ \hline
$\GR{E}{7}$  & $\vp_1$  & yes & $(1+q^5)(1+q^9)(1+q+\dots +q^{13})$ \\ \hline
$\GR{G}{2}$  & $\vp_1$  & no & $1+q+\dots +q^6$
\\ \hline
\end{tabular}
\end{center}
\end{table}
One can observe that each $\lb$ in Table~1 is a multiple of
a fundamental weight.
Using polynomials $\mlmq$, one can give a conceptual proof of this fact.
This will appear elsewhere.

\section{A connection with (equivariant) cohomology}
\setcounter{equation}{0}
Let $\ov{\cc}$ be an Artinian graded commutative associative $\bbk$-algebra,
$\ov{\cc}=\oplus_{i=0}^d\ov{\cc}_i$.
Suppose $\dim\ov{\cc}_d=1$, and let $\xi$ be a nonzero linear form
on $\ov{\cc}$ that annihilates the space
$\ov{\cc}_0\oplus\dots\oplus\ov{\cc}_{d-1}$. Then $\ov{\cc}$ is
called a {\it Poincar\'e duality algebra}, if the bilinear form
$(x,y)\mapsto \xi(xy)$, $x,y\in\ov{\cc}$, is non-degenerate.
This name suggests that $\ov{\cc}$ looks very much as if it were
the cohomology algebra of some `good' manifold.
It is easily seen that $\ov{\cc}$ is a Poincar\'e duality algebra if
and only if $\dim{\mathrm {soc}}(\ov{\cc})=1$, i.e., $\ov{\cc}$ is
Gorenstein.
\\[.6ex]
In the previous section, we have proved that
$\cR^{(\lb)}=
(\evl)^A$ is an Artinian  Gorenstein $\bbk$-algebra, if $\vlb$ is
{\sf wmf}. (Here $d=2\,\hot(\lb)$, which is not necessarily even.)
In this situation, one has more evidences in favour of the assertion
that $\cR^{(\lb)}$ can be
a cohomology algebra. Recall that $e$ is a sum of root vectors corresponding
to the simple roots. Hence $e(v_{-\lb^*})\in (\vlb)_1$, i.e., $e\in\cR^{(\lb)}_1$.

\begin{s}{Proposition}
The multiplication operator $e: \cR^{(\lb)}_i\to \cR^{(\lb)}_{i+1}$ is
injective for $i \le [(2\hot(\lb)-1)/2]$ and surjective for
$i\ge [\hot(\lb)]$.
\end{s}\begin{proof}
This follows from the $\tri$-theory and the equality $\dlq(q)=\cf_\lb(q)$.
\end{proof}%
Hence if $\cR^{(\lb)}=H^*(X)$ for some algebraic
variety $X$, then $e\in \cR^{(\lb)}_1$
can be regarded as the class of a hyperplane section, and the hard Lefschetz
theorem holds for $X$.
If $G$ is simple, then $\lb$ is necessarily fundamental. Therefore
the subspace $(\vlb)_1\subset\vlb$ is 1-dimensional. Indeed, if
$\ap\in\Pi$ is the unique root such that $(\ap,\lb^*)\ne 0$, then
$(\vlb)_1=\vlb^{\ap-\lb^*}$. In terms of $\cR^{(\lb)}$ this means
that $\cR^{(\lb)}_1=\bbk e$.
%
That is, such $X$ should satisfy the condition $b_2(X)=1$.
\\[.7ex]
At the rest of the paper, we assume that $\bbk=\Bbb C$, and
consider cohomology
with complex coefficients.
Now we give a hypothetical description of such $X$ in case $G$ is
simple. More precisely, we conjecture that there exists
$X_\lb\subset {\Bbb P}(\vlb)$ such that odd cohomology of $X_\lb$ vanishes
and  $\cR^{(\lb)}\simeq H^*(X_\lb)$; moreover,
the $\g$-endomorphism algebra
$\clg$ gives the equivariant cohomology of $X_\lb$.
We refer to \cite{brion} for a nice introduction to equivariant
cohomology. As usual, $\vlb$ is a {\sf wmf} $G$-module.
The variety we are seeking for should satisfy the constraints
$\chi(X_\lb)=\dim\vlb$ and $\dim X_\lb= 2\,\hot(\lb)$.
Notice that the set of $T$-fixed points in ${\Bbb P}(\vlb)$ is finite:
\[
{\Bbb P}(\vlb)^T=\bigcup_{\nu\dashv\vlb} \langle{v^\nu}\rangle \ ,
\]
where $\langle{v^\nu}\rangle$ is the image of $v^\nu\in\vlb^\nu$ in the
projective space.
Since $\chi(Y)=\chi(Y^T)$ for any algebraic variety $Y$ acted upon by a
torus $T$ \cite{euler}, we see that our variety $X_\lb$ must contain all
$\langle{v^\nu}\rangle$, $\nu\dashv\vlb$.
This provides some explanation for the following description.
Let $\mu\dashv\vlb$ be the unique dominant minuscule weight.
For instance, $\mu=0$ if and only if $\lb\in\cq$.

{\bf -- }
If $(G,\lb)\ne (\GR{A}{n}, m\vp_1)$ or $(\GR{A}{n}, m\vp_n)$ with $m>n+1$, then
we define $X_\lb$ to be the closure of the $G$-orbit of the line
$\langle{v_\mu}\rangle$ in the projective space ${\Bbb P}(\vlb)$.

{\bf -- }
Alternatively, for $G=\GR{A}{n}$ and $\lb=m\vp_1$ with any $m\ge 1$,
let $X_\lb$ be the projectivisation of the
variety of decomposable forms of degree $m$ in
$\VV_{m\vp_1}=S^m (\VV_{\vp_1})$, where a form of degree $m$ is said to be
{\it decomposable}, if it is a product of $m$ linear forms.
(Similar definition applies to $\lb=m\vp_n$.)
\\
As is easily seen,  the two constructions coincide, if they
both apply, i.e. for $(\GR{A}{n}, m\vp_1 {\textrm { or }} m\vp_n)$ with
$m\le n+1$. But I do not see how to give a uniform description of $X_\lb$
in all cases. Direct calculations show that $\dim X_\lb=2\,\hot(\lb)$;
e.g. for $(\GR{A}{n}, m\vp_1)$ we have $\dim X_\lb=mn$.

\begin{s}{Conjecture} Let $\vlb$ be a {\sf wmf} $G$-module. Then \\
(1) the variety $X_\lb\subset
{\Bbb P}(\vlb)$ is rationally smooth,
odd cohomology of $X_\lb$ vanishes, and
$H^*(X_\lb)\simeq \cR^{(\lb)}$. In particular,
the Poincar\'e polynomial of $H^*(X_\lb)$ is equal to $\dlq(q^2)$.
\\
(2) Let $G_c\subset G$ be a maximal compact subgroup of $G$. Then
$H^*_{G_c}(X_\lb)\simeq \clg$.
\end{s}%
Let $J_+\subset J$ be the augmentation ideal.
Since  $\clg/J_+\clg\simeq \cR^{(\lb)}$, $J\simeq H^*_{G_c}(\{pt\})$,
and $H^*_{G_c}(X_\lb)/(J_+)\simeq H^*(X_\lb)$ (see \cite[Prop.\,2]{brion}),
the first part of the conjecture follows from the second one.
Actually, the conjecture is true for most of items in Table 1, in
particular, for all minuscule weights.
We list below all known results supporting the conjecture.

$\bullet$ \ If $\lb$ is minuscule, then $\lb=\mu$ and therefore
$X_\lb=G/P_\lb$, a generalised flag
variety. Here the conjecture follows from Theorem~\ref{main2}
and the well-known description of $H^*_{G_c}(G/P_\lb)$.
Indeed, if $P_\lb$ is an {\it arbitrary\/} parabolic subgroup (i.e. $\lb$
is not necessarily minuscule), then $H^*_{G_c}(G/P_\lb)
\simeq \bbk[\te]^{W_\lb}$.

$\bullet$ \ For the simplest representations of $\GR{B}{n}$ and $\GR{G}{2}$,
we have $\mu=0$ and $X_\lb= {\Bbb P}(\vlb)$. On the other hand, the
formulae for $\dlq(q)$ in Table~1 shows that
$(\evl)^A$ is generated by $e$ as $\bbk$-algebra, and
the equality $H^*(X_\lb)= \cR^{(\lb)}$  follows. This also shows that
$\clg$ and $H^*_{G_c}(X_\lb)$ have the same Poincar\'e series.

$\bullet$ \  We have essentially four cases with non-minuscule $\lb$:
$(\GR{B}{n}, \vp_1),\,(\GR{C}{3}, \vp_3),\,
(\GR{G}{2}, \vp_1)$, and $(\GR{A}{n}, m\vp_1)$, $m\ge 2$.
For the first three cases and for the last one with $m=2$,
$X_\lb$ is a compact multiplicity-free $G_c$-space; in other words,
$X_\lb$ is a spherical $G$-variety.
Therefore making use
of Theorem~9 in \cite{brion}, one obtains the description
of $H^*_{G_c}(X_\lb)$. In these ``spherical'' cases, the number of
dominant weights in $\vlb$ equals 2. Therefore the structure of $\clg$
is not too complicated, and
this can be used for proving that
$\clg=H^*_{G_c}(X_\lb)$.

$\bullet$ \ The projectivisation of the variety of decomposable
forms of degree $m$ in $n+1$ variables is isomorphic with
$({\Bbb P}^{n})^m/\Sigma_m$, where $\Sigma_m$ is the symmetric group
permuting factors in ${\Bbb P}^{n}\times\dots\times {\Bbb P}^{n}$
($m$ times), see \cite[Theorem\,1.3]{brion2}, \cite[4.2]{triple}.
Notice that
\[
H^*(({\Bbb P}^n)^m/\Sigma_m)\simeq H^*(({\Bbb P}^n)^m)^{\Sigma_m}
\simeq
(\bbk[x_1,\dots,x_m]/(x_1^{n+1},
\dots x_m^{n+1}))^{\Sigma_m} \ ,
\]
the algebra of truncated symmetric polynomials. It is easily seen that
its dimension is equal to
$\genfrac{(}{)}{0pt}{}{m+n}{m}=\dim \VV_{m\vp_1}$. A somewhat more
bulky but still elementary calculation shows that the Poincar\'e
polynomial of $H^*(({\Bbb P}^n)^m/\Sigma_m)$ is equal to the Dynkin
polynomial for $(\GR{A}{n}, m\vp_1)$. Our proof is purely
combinatorial. We establish a natural one-to-one correspondence between
suitably chosen bases of $H^*(({\Bbb P}^n)^m/\Sigma_m)$ and $\VV_{m\vp_1}$
such that $H^i(({\Bbb P}^n)^m/\Sigma_m)$ corresponds to
$(\VV_{m\vp_1})_i$.
It is not, however, clear how
to compare the multiplicative structure of $H^*(({\Bbb P}^n)^m/\Sigma_m)$
and $(\End\VV_{m\vp_1})^A$.



\begin{thebibliography}{LSo1}

\bibitem[BiB]{euler} {\sc A.~Bialynicki-Birula}. On fixed points
schemes of actions of multiplicative and additive groups,
{\it Topology}, {\bf 12}(1973), 99--103.

\bibitem[Bri1]{brion2} {\sc M.~Brion}. Stable properties of plethysm:
on two conjectures of Foulkes. {\it Manuscripta Math.}
{\bf 80}(1993), 347--371.

\bibitem[Bri2]{brion} {\sc M.~Brion}. Equivariant cohomology and
equivariant intersection theory, In: ``Representation Theories and
Algebraic Geometry'', A.~Broer (Ed.), 1--37. Kluwer, Dordrecht/Boston/London,
1998.

\bibitem[Bro1]{bram93} {\sc A.~Broer}. Line bundles on the
cotangent bundle of the flag variety, {\it Invent. Math.}
{\bf 113}(1993), 1--20.

\bibitem[Bro2]{bram} {\sc A.~Broer}. The sum of generalized exponents and
Chevalley's restriction theorem, {\it Indag. Math.} {\bf 6}(1995), 385--396.

\bibitem[Bry]{jump} {\sc R.K.~Brylinski}.  Limits of weight spaces, Lusztig's
$q$-analogs, and fiberings of adjoint orbits,
{\it J. Amer. Math. Soc.} {\bf 2}(1989), 517--533.

\bibitem[GKZ]{triple} {\sc I.M.~Gelfand, M.~Kapranov, A.~Zelevinsky}.
``Discriminants, Resultants, and Multidimensional Determinants'',
Birkh\"auser, Boston/Basel/Berlin, 1994.

\bibitem[Gu1]{ranee87} {\sc R.K.~Gupta}. Characters and $q$-analog of weight
multiplicity, {\it Bull. London Math. Soc.}(2) {\bf 36}(1987), 68--76.

\bibitem[Gu2]{ranee} {\sc R.K.~Gupta}. Generalized exponents via
Hall-Littlewood symmetric functions, {\it Bull. Amer. Math. Soc.}
{\bf 16}(1987), 287--291.

\bibitem[Dy]{vereteno} {\sc E.B.~Dynkin}.
Some properties of the weight system  of a linear representation
of a semisimple Lie group, {\it Doklady Akad. Nauk SSSR}, {\bf 71}(1950),
{\rus N0}\,2, 221--224 (Russian).


\bibitem[Gr]{graham} {\sc W.~Graham}. Functions on the universal cover of
the principal nilpotent orbit, {\it Invent. Math.} {\bf 108}(1992), 15--27.

\bibitem[He]{wim} {\sc W.~Hesselink}. Characters of the nullcone,
{\it Math. Annalen}, {\bf 252}(1980), 179--182.

\bibitem[Ho]{howe}
{\sc R.~Howe}.  Perspectives on Invariant Theory: Schur duality,
multiplicity-free actions and beyond,
In:  ``The Schur lectures (1992)'',
Israel Math. Conf. Proceedings, {\bf 8}(1995), 1--182.

\bibitem[Hum]{hump}  {\sc J.E.~Humphreys}. ``Reflection Groups and
Coxeter Groups'', Cambridge Univ. Press, 1992.

\bibitem[Ki]{family} {\sc A.A.~Kirillov}. Introduction to family algebras,
{\it Moscow Math. J.} {\bf 1}(2001), 49--64.

\bibitem[Ko]{ko63} {\sc B.~Kostant}. Lie group representations in
polynomial rings, {\it  Amer. J. Math.} {\bf 85}(1963), 327--404.

\bibitem[Pa]{ya} {\sc D.~Panyushev}.
On covariants of reductive algebraic groups, {\it Indag. Math.},
to appear.

\bibitem[PRV]{prv} {\sc K.R.~Parthasarathy, R.~Ranga Rao, V.S.~Varadarajan}.
Representations of complex semisimple Lie groups and Lie algebras,
{\it Annals Math.} {\bf 85}(1967), 383--429.

\bibitem[Sta1]{unimod} {\sc R.P.~Stanley}. Unimodal sequences arising from
Lie algebras, In: Young Day Proceedings, T.V.~Narayana et al. (Eds.), 127--136.
Dekker, New York/Basel, 1980.

\bibitem[Sta2]{combin} {\sc R.P.~Stanley}. ``Combinatorics and commutative
algebra'' (Second Edition), Boston/Basel/Berlin, Birkh\"auser, 1996.

\bibitem[Ste]{stembridge}  {\sc J.~Stembridge}. On minuscule representations,
plane partitions and involutions in complex Lie groups,
{\it Duke math. J.} {\bf 73}(1994), 469--490.

\bibitem[VO]{vion}
{\sc E.B.\,Vinberg} and {\sc A.L.\,Onishchik}.
``Seminar on Lie groups and
algebraic groups",  Moskva: ``Nauka'' 1988 (Russian).
English translation: {\sc A.L.\,Onishchik} and {\sc E.B.\,Vinberg}:
``Lie groups and algebraic groups''. Berlin Heidelberg
New York: Springer 1990.

\end{thebibliography}
\end{document}